# Random walk models associated with distributed fractional order differential equations

Sabir Umarov[1,*,†] and Stanly Steinberg[1]

*University of New Mexico*

**Abstract:** In this paper the multi-dimensional random walk models governed by distributed fractional order differential equations and multi-term fractional order differential equations are constructed. The scaling limits of these random walks to a diffusion process in the sense of distributions is proved.

## 1. Introduction

In this paper we construct new random walks connected with fractional order differential equations. Namely, the governing equations corresponding to the constructed random walks are multi-term or distributed fractional order differential equations. Nowadays the connection between random walk and fractional order dynamics is well known, see, for instance [1, 17, 26, 33, 38]. A number of constructive random walk models governed by fractional differential equations in the one-dimensional case were studied by Gillis, et al. [12], Chechkin, et al. [7], Gorenflo, et al. [15, 16], and in the $n$-dimensional case by Umarov [35], Umarov, et al. [36], Andries, et al. [3]. The governing equation in these studies depends on parameters $\beta \in (0,1]$ and $\alpha \in (0,2]$, and is given by the fractional order differential equation

$$(1) \qquad \mathcal{D}^\beta u(t,x) = D_0^\alpha u(t,x),\ t > 0,\ x \in \mathbb{R}^N,$$

where $\mathcal{D}^\beta$ is the time-fractional derivative in some sense, and $D_0^\alpha, 0 < \alpha < 2$, is the pseudo-differential operator with the symbol $-|\xi|^\alpha$. The precise definitions will be given below.

In the present paper we construct the random walks the governing equation of which is a distributed space fractional order differential equation

$$(2) \qquad \frac{\partial}{\partial t} u(t,x) = \int_0^2 a(\alpha) D_0^\alpha u(t,x) d\alpha,\ t > 0,\ x \in \mathbb{R}^N,$$

where $a(\alpha)$ is a positive integrable function (positively defined distribution).

The study of properties of distributed order differential operators and their applications to diffusion processes has been developed extensively in recent years, although such operators were first mentioned by Caputo [5, 6] in 1960th. The distributed order differential equations have been used by Bagley, et al. [4] to model

[1]Mathematics and Statistics Department, University of New Mexico, Albuquerque, NM, e-mail: sabir@math.unm.edu; stanly@math.unm.edu

*Corresponding author.

†Supported by Fulbright Program.

the input-output relationship of linear time-variant system, by Lorenzo, et al. [22] to study of rheological properties of composite materials, by Chechkin, et al. [8] to model some ultraslow and lateral diffusion processes. Diethelm, et al. [9] studied the numerical aspects of such equations. Umarov, et al. [37] studied general properties of distributed order differential equations and solvability problems of the Cauchy and multipoint value problems.

The method used in this paper for construction of multi-dimensional random walks are essentially based on the symbolic calculus of pseudo-differential operators and on the convergence properties of some simple cubature formulas. This method is new even in the one-dimensional case and was suggested recently in [35, 36]. We note that the scaling limits are obtained in terms of characteristic functions of transition probabilities. The equivalence of corresponding convergence notions is well-known (see, [11, 13]). See also the recent book by M. Meerschaert and Scheffler where the multi-dimensional operator stable probability distributions are studied and analogs of different type limits considered. Multi-dimensional random walks are frequently used in modeling various processes in different areas [1, 2, 24, 25, 26, 31].

The present report is organized as follows. In Section 2 we give preliminaries simultaneously introducing the terminology that will be used in the paper. We also recall some properties of pseudo-differential operators with constant symbols and lay out some elementary properties of symbols. These properties play an essential role later in the study of the diffusion limits of random walks. In Section 3 we formulate our random walk problem in terms of sequences of i.i.d. (independent identically distributed) random vectors. In this Section we also formulate the main results.

## 2. Preliminaries

We use the following notation. $\mathbb{R}^N$ is the $N$-dimensional Euclidean space with coordinates $x = (x_1, \ldots, x_N)$; $\mathbb{Z}^N$ is the $N$-dimensional integer-valued lattice with nodes $j = (j_1, \ldots, j_N)$. We denote by $x_j = (hj_1, \ldots, hj_N), j \in \mathbb{Z}^N$, the nodes of the uniform lattice $\mathbb{Z}^N_h$ defined as $(h\mathbb{Z})^N$ with a positive number $h$, the mesh width.

We assume that a walker is located at the origin $x_0 = (0, \ldots, 0)$ at the initial time $t = 0$. In our random walk, at every time instants $t_1 = \tau, t_2 = 2\tau, \ldots, t_n = n\tau, \ldots$ the walker jumps through the nodes of the lattice $\mathbb{Z}^N_h$. By $p_j, j \in \mathbb{Z}^N$, we denote transition probabilities. Namely, $p_j$ means a probability of jumping of the walker from a point $x_k \in \mathbb{Z}^N_h$ to a point $x_{j+k} \in \mathbb{Z}^N_h$, where $j$ and $k$ are in $\mathbb{Z}^N$. Transition probabilities satisfy the non-negativity and normalization conditions:

1. $p_j \geq 0, j \in \mathbb{Z}^N$;
2. $\sum_{j \in \mathbb{Z}^N} p_j = 1$.

Transition probabilities $\{p_j, j \in \mathbb{Z}^N\}$ are associated with a discrete function $p : \mathbb{Z}^N \to [0, 1]$. For given two transition probabilities, $p$ and $q$ we define the convolution operation $p * q$ by

$$(p * q)_j = \sum_{k \in \mathbb{Z}^N} p_k q_{j-k}, j \in \mathbb{Z}^N.$$

Let $f$ be a continuous function integrable over $\mathbb{R}^N$. Then, as is known [32], the rectangular cubature formula

$$(3) \qquad \int_{\mathbb{R}^N} f(x)dx = h^N \sum_{j \in \mathbb{Z}^N} f(x_j) + o(1)$$

is valid.

The operators in our random walk models have a close relationship to pseudo-differential operators with the symbols depending only on the dual variables. Symbols are allowed to have singularities. For general orientation to the theory of such operators we refer to [10, 14, 18, 19, 20, 21, 34].

Let $A(D), D = (D_1, \ldots, D_N), D_j = \partial/i\partial x_j, j = 1, \ldots, N$, be a pseudo-differential operator with a symbol $A(\xi)$ not depending on $x$, and defined in $\mathbb{R}^N$. We refer to the variable $\xi$ as a dual variable. Both type of variables, $x$ and $\xi$ belong to $\mathbb{R}^N$ (more precisely, $\xi$ belongs to the conjugate $(\mathbb{R}^N)^* = \mathbb{R}^N$). To avoid confusion sometimes we write $\mathbb{R}^N_x$ and $\mathbb{R}^N_\xi$, indicating their relevance to the variables $x$ and $\xi$ respectively. Further, for a test function $\varphi(x)$ taken from the classical space $S(\mathbb{R}^N_x)$ of rapidly decreasing functions, the Fourier transform

$$\hat{\varphi}(\xi) = F[\varphi](\xi) = \int_{\mathbf{R}^N} \varphi(x) e^{-ix\xi} dx$$

is well defined and belongs again to $S(\mathbb{R}^N_\xi)$. Let $S'(\mathbb{R}^N)$ be the space of tempered distributions, i.e. the dual space to $S(\mathbb{R}^N)$. The Fourier transform for distributions $f \in S'(\mathbb{R}^N_x)$ is usually defined by the extension formula $(\hat{f}(\xi), \varphi(\xi)) = (f(x), \hat{\varphi}(x))$, with the duality pairing $(.,.)$ of $S'(\mathbb{R}^N_\xi)$ and $S(\mathbb{R}^N_\xi)$.

**Definition 2.1.** Assume $G$ to be an open domain in $\mathbb{R}^N_\xi$. Let a function $f$ be continuous and bounded on $\mathbb{R}^N_x$ and have a Fourier transform (taken in the sense of distributions) $\hat{f}(\xi)$ with compact support in $G$. We denote by $\Psi_G(\mathbb{R}^N_x)$ the set of all such functions endowed with the following convergence. A sequence of functions $f_m \in \Psi_G(\mathbb{R}^N_x)$ is said to converge to an element $f_0 \in \Psi_G(\mathbb{R}_x^N)$ iff:

1. there exists a compact set $K \subset G$ such that $\operatorname{supp} \hat{f}_m \subset K$ for all $m = 1, 2, \ldots$;
2. $\|f_m - f_0\| = \sup |f_m - f_0| \to 0$ for $m \to \infty$.

In the case $G = \mathbb{R}^N_\xi$ we write simply $\Psi(\mathbb{R}_x{}^N)$ omitting $\mathbb{R}^N_\xi$ in the index of $\Psi_G(\mathbb{R}^N_x)$.

Note that according to the Paley-Wiener theorem functions in $\Psi_G(\mathbb{R}^N_x)$ are entire functions of finite exponential type (see [27], [10]).

Denote by $H^s(\mathbb{R}^N_x), s \in (-\infty, +\infty)$ the Sobolev space of elements $f \in S'(\mathbb{R}^N_x)$ for which $(1 + |\xi|^2)^{s/2}|\hat{f}(\xi)| \in L_2(\mathbb{R}^N_\xi)$. It is known [18] that if $f \in L_p(\mathbb{R}_x{}^N)$ with $p > 2$, then its Fourier transform $\hat{f}$ belongs to $H^{-s}(\mathbb{R}^N_\xi), s > N(\frac{1}{2} - \frac{1}{p})$. Letting $p \to \infty$ we get $\hat{f} \in H^{-s}(\mathbb{R}^N_\xi), s > \frac{N}{2}$ for $f \in L_\infty(\mathbb{R}_x{}^N)$. It follows from this fact and the Paley-Wiener theorem that the Fourier transform of $f \in \Psi_G(\mathbb{R}^N)$ belongs to the space

$$\bigcap_{s > \frac{N}{2}} H_c^{-s}(G),$$

where $H_c^{-s}(G)$ is a negative order Sobolev space of functionals with compact support on $G$. Hence $\hat{f}$ is a distribution, which is well defined on continuous functions.

Let $\Psi'_{-G}(\mathbb{R}^N)$ be the space of all linear bounded functionals defined on the space $\Psi_G(\mathbb{R}^N)$ endowed with the weak (dual with respect to $\Psi_G(\mathbb{R}^N)$) topology. By the weak topology we mean that a sequence of functionals $g_m \in \Psi'_{-G}(\mathbb{R}^N)$ converges to an element $g_0 \in \Psi'_{-G}(\mathbb{R}^N)$ in the weak sense if for all $f \in \Psi_G(\mathbb{R}^N)$ the sequence of numbers $\langle g_m, f \rangle$ converges to $\langle g_0, f \rangle$ as $m \to \infty$. By $\langle g, f \rangle$ we mean the value of $g \in \Psi'_{-G}(\mathbb{R}^N)$ on an element $f \in \Psi_G(\mathbb{R}^N)$.

**Definition 2.2.** Let $A(\xi)$ be a continuous function defined in $G \subset \mathbb{R}^N_\xi$. A pseudo-differential operator $A(D)$ with the symbol $A(\xi)$ is defined by the formula

$$(4) \qquad A(D)\varphi(x) = \frac{1}{(2\pi)^N}(\hat{\varphi}, A(\xi)e^{-ix\xi}), \ \varphi \in \Psi_G(\mathbb{R}^N).$$

Obviously, the function $A(\xi)e^{-ix\xi}$ is continuous in $G$. Thus, $A(D)$ in Eq. (4) is well defined on $\Psi_G(\mathbb{R}^N)$. If $\hat{\varphi}$ is an integrable function with $\operatorname{supp} \hat{\varphi} \subset G$, then (4) takes the usual form of pseudo-differential operator

$$A(D)\varphi(x) = \frac{1}{(2\pi)^N} \int A(\xi)\hat{\varphi}(\xi)e^{-ix\xi}d\xi,$$

with the integral taken over $G$. Note that in general this integral may not make sense even for infinitely differentiable functions with finite support (see [14]).

Now we define the operator $A(-D)$ acting in the space $\Psi'_{-G}(\mathbb{R}^N)$ by the extension formula

$$(5) \qquad < A(-D)f, \varphi > = < f, A(D)\varphi >, \ f \in \Psi'_{-G}(\mathbb{R}^N), \ \varphi \in \Psi_G(\mathbb{R}^N).$$

We recall (see [14]) some basic properties of pseudo-differential operators introduced above.

**Lemma 2.3.** *The pseudo-differential operators $A(D)$ and $A(-D)$ with a continuous symbol $A(\xi)$ act as*

1. $A(D) : \Psi_G(\mathbb{R}^N) \to \Psi_G(\mathbb{R}^N)$,
2. $A(-D) : \Psi'_{-G}(\mathbb{R}^N) \to \Psi'_{-G}(\mathbb{R}^N)$

*respectively, and are continuous.*

**Lemma 2.4.** *Let $A(\xi)$ be a function continuous on $\mathbb{R}^N$. Then for $\xi \in \mathbb{R}^N$*

$$A(D)\{e^{-ix\xi}\} = A(\xi)e^{-ix\xi}.$$

*Proof.* For any fixed $\xi \in \mathbb{R}^N$ the function $e^{-ix\xi}$ is in $\Psi(\mathbb{R}^N)$. We have

$$A(D)\{e^{-ix\xi}\} = \frac{1}{(2\pi)^N} \int_{\mathbb{R}^N} A(\eta) e^{-ix\eta} d\mu_\xi(\eta),$$

where $d\mu_\xi(\eta) = F_\eta[e^{-ix\xi}]d\eta = (2\pi)^N \delta(\eta - \xi)d\eta$. Hence $A(D)\{e^{-ix\xi}\} = A(\xi)e^{-ix\xi}$. $\square$

**Corollary 2.5.**
1. $A(\xi) = (A(D)e^{-ix\xi})e^{ix\xi}$;
2. $A(\xi) = (A(D)e^{-ix\xi})_{|x=0}$;
3. $A(\xi) = < A(-D)\delta(x), e^{-ix\xi} >$, where $\delta$ is the Dirac distribution.

**Remark 2.6.** Since the function $e^{-ix\xi}$ does not belong to $S(\mathbb{R}^N)$ and $D(\mathbb{R}^N)$, the representations for the symbol obtained in Lemma 2.4 and Corollary 2.5 are not directly applicable in these spaces.

Denote by $D_0^\alpha$, $0 < \alpha \leq 2$, the pseudo-differential operator with the symbol $-|\xi|^\alpha$. It is evident that $D_0^\alpha$ coincides with the Laplace operator $\Delta$ for $\alpha = 2$. For $\alpha < 2$ it can be considered as a fractional power of the Laplace operator, namely $D_0^\alpha = -(-\Delta)^{\alpha/2}$. $D_0^\alpha$ can also be represented as a hypersingular integral (see, e.g. [29])

$$(6) \qquad D_0^\alpha f(x) = b(\alpha) \int_{\mathbb{R}^N_y} \frac{\Delta_y^2 f(x)}{|y|^{N+\alpha}} dy,$$

where $\Delta_y^2$ is the second order centered finite difference in the $y$ direction, and $b(\alpha)$ is norming constant defined as

$$(7) \qquad b(\alpha) = \frac{\alpha \Gamma(\frac{\alpha}{2}) \Gamma(\frac{N+\alpha}{2}) \sin \frac{\alpha \pi}{2}}{2^{2-\alpha} \pi^{1+N/2}}.$$

It is seen from (7) that in the representation (6) of $D_0^\alpha$ the value $\alpha = 2$ is singular.

**Lemma 2.7.** *For the symbol of $D_0^\alpha$ the following equalities hold true:*

$$(8) \qquad \left(D_0^\alpha e^{ix\xi}\right)|_{x=0} = b(\alpha) \int_{\mathbb{R}_y^N} \frac{\Delta_y^2 e^{ix\xi}}{|y|^{N+\alpha}} dy|_{x=0} = -|\xi|^\alpha, 0 < \alpha < 2.$$

*Proof.* This statement is a direct implication of Corollary 2.5 applied to the operator $D_0^\alpha$. □

The cubature formula (3) yields for the integral in the right hand side of (6)

$$(9) \qquad \int_{\mathbb{R}^N} \frac{\Delta_y^2 f(x_j)}{|y|^{N+\alpha}} dy = h^\alpha \sum_{k \in \mathbb{Z}^N} \frac{\Delta_k^2 f_j}{|k|^{N+\alpha}} + o(1), \ j \in \mathbb{Z}^N,$$

where $f_j = f(x_j)$ and $|k|$ is Euclidean norm of $k = (k_1, \ldots, k_N) \in \mathbb{Z}^N$.

Consider the distributed space fractional order differential equation

$$(10) \qquad \frac{\partial}{\partial t} u(t,x) = \int_0^2 a(\alpha) D_0^\alpha u(t,x) d\alpha, \ t > 0, \ x \in \mathbb{R}^N,$$

where $a(\alpha)$ is a positive (in general, generalized) function defined in $(0, 2]$. A distribution $G(t, x)$, which satisfies the equation (10) in the weak sense and the condition

$$(11) \qquad G(0, x) = \delta(x), x \in \mathbb{R}^N,$$

where $\delta(x)$ is the Dirac's distribution, is called a fundamental solution of the Cauchy problem (10), (11). In the particular case of

$$(12) \qquad a(\alpha) = \sum_{m=1}^M a_m \delta(\alpha - \alpha_m), \ 0 < \alpha_1 < \cdots < \alpha_M \leq 2,$$

with positive constants $a_m$ we get a multiterm space fractional differential equation

$$(13) \qquad \frac{\partial}{\partial t} u(t,x) = \sum_{m=1}^M a_m D_0^{\alpha_m} u(t,x) \, t > 0, \ x \in \mathbb{R}^N.$$

Denote the operator on the right hand side of the equation (10) by $\mathcal{B}(D)$. It can be represented as a pseudo-differential operator with the symbol

$$(14) \qquad \mathcal{B}(\xi) = -\int_0^2 a(\alpha) |\xi|^\alpha d\alpha.$$

It is not hard to verify that the fundamental solution of equation (10) is

$$(15) \qquad G(t,x) = F^{-1}\left(e^{t\mathcal{B}(\xi)}\right),$$

where $F^{-1}$ stands for the inverse Fourier transform. In the particular case of $a(\alpha) = \delta(\alpha - 2)$ we have the classical heat conduction equation

$$\frac{\partial}{\partial t} u(t,x) = \Delta u(t,x),\ t > 0,\ x \in \mathbf{R}^N,$$

whose fundamental solution is the Gaussian probability density function evolving in time

$$G_2(t,x) = \frac{1}{(4\pi t)^{n/2}} e^{\frac{-|x|^2}{4t}}.$$

For $a(\alpha) = \delta(\alpha - \alpha_0)$, $0 < \alpha_0 < 2$, the corresponding fundamental solution is the Levy $\alpha_0$-stable probability density function [30]

(16) $$G_{\alpha_0}(t,x) = \frac{1}{(2\pi)^N} \int_{\mathbb{R}^N} e^{-t|\xi|^{\alpha_0}} e^{ix\xi} d\xi.$$

The power series representation of the stable Levy probability density functions is studied in [3, 23, 33]. Recall also that $\alpha_0 = 1$ corresponds to the Cauchy-Poisson probability density (see [28])

$$G_1(t,x) = \frac{\Gamma(\frac{n+1}{2})}{\pi^{(n+1)/2}} \frac{1}{(|x|^2 + t^2)^{(n+1)/2}}.$$

## 3. Main results: construction of random walks

In this Section we construct random walks associated with distributed space fractional order differential equations (10). More precisely, we show that the special scaling limit of the constructed random walk is a diffusion process whose probability density function is the fundamental solution of (10).

Let $\mathbf{X}$ be an N-dimensional random vector [25] which takes values in $\mathbb{Z}^N$. Let the random vectors $\mathbf{X}_1, \mathbf{X}_2, \ldots$ also be N-dimensional independent identically distributed random vectors, all having the same probability distribution, common with $\mathbf{X}$. We introduce a spatial grid $\{x_j = jh, j \in \mathbb{Z}^N\}$, with $h > 0$ and temporal grid $\{t_n = n\tau, n = 0, 1, 2, \ldots\}$ with a step $\tau > 0$. Consider the sequence of random vectors

$$\mathbf{S}_n = h\mathbf{X}_1 + h\mathbf{X}_2 + \cdots + h\mathbf{X}_n,\ n = 1, 2, \ldots$$

taking $\mathbf{S}_0 = \mathbf{0}$ for convenience. We interpret $\mathbf{X}_1, \mathbf{X}_2, \ldots,$ as the jumps of a particle sitting in $x = x_0 = \mathbf{0}$ at the starting time $t = t_0 = 0$ and making a jump $\mathbf{X}_n$ from $\mathbf{S}_{n-1}$ to $\mathbf{S}_n$ at the time instant $t = t_n$. Then the position $\mathbf{S}(t)$ of the particle at time $t$ is

$$\sum_{1 \leq k \leq t/\tau} \mathbf{X}_k.$$

Denote by $y_j(t_n)$ the probability of sojourn of the walker at $x_j$ at the time $t_n$. Taking into account the recursion $\mathbf{S}_{n+1} = \mathbf{S}_n + h\mathbf{X}_n$ we have

(17) $$y_j(t_{n+1}) = \sum_{k \in \mathbb{Z}^N} p_k y_{j-k}(t_n),\ j \in \mathbb{Z}^N,\ n = 0, 1, \ldots$$

The convergence of the sequence $\mathbf{S}_n$ when $n \to \infty$ means convergence of the discrete probability law $(y_j(t_n))_{j \in \mathbb{Z}^N}$, properly rescaled as explained below, to the probability law with a density $u(t,x)$ in the sense of distributions (in law). This is equivalent to the locally uniform convergence of the corresponding characteristic functions (see for details [25]). We use this idea to prove the convergence of the sequence of characteristic functions of the constructed random walks to the fundamental solution of distributed order diffusion equations.

In order to construct a random walk relevant to (10) we use the approximation (3) for the integral on the right hand side of (6), namely

$$D_0^\alpha u(t, x_j) \approx b(\alpha) \sum_{k \in \mathbb{Z}^N} \frac{u_{j+k}(t) - 2u_j(t) + u_{j-k}(t)}{|k|^{N+\alpha} h^\alpha},$$

and the first order difference ratio

$$\frac{\partial u}{\partial t} \approx \frac{u_j(t_{n+1}) - u_j(t_n)}{\tau}$$

for $\frac{\partial u}{\partial t}$ with the time step $\tau = t/n$. Then from (10) we derive the relation (17) with the transition probabilities

$$(18) \qquad p_k = \begin{cases} 1 - 2\tau \sum_{m \neq 0} \frac{Q_m(h)}{|m|^N}, & \text{if } k = 0; \\ 2\tau \frac{Q_m(h)}{|k|^N}, & \text{if } k \neq 0, \end{cases}$$

where

$$(19) \qquad Q_m(h) = \int_0^2 |m|^{-\alpha} \rho(\alpha, h) d\alpha, \quad \rho(\alpha, h) = \frac{a(\alpha) b(\alpha)}{h^\alpha}.$$

Assume that the condition

$$(20) \qquad \sigma(\tau, h) := 2\tau \sum_{m \neq 0} \frac{Q_m(h)}{|m|^N} \leq 1.$$

is fulfilled. Then, obviously, the transition probabilities satisfy the properties:

1. $\sum_{k \in \mathbb{Z}^N} p_k = 1$;
2. $p_k \geq 0, k \in \mathbb{Z}^N$.

Introduce the function

$$\mathcal{R}(\alpha) = \sum_{k \neq 0} \frac{1}{|k|^{N+\alpha}} = \sum_{m=1}^\infty \frac{M_m}{m^{N+\alpha}}, \quad 0 < \alpha \leq 2,$$

where $M_m = \sum_{|k|=m} 1$. (In the one-dimensional case $\mathcal{R}(\alpha)$ coincides with the Riemann's zeta-function, $\mathcal{R}(\alpha) = 2\zeta(1+\alpha)$.) The Eq. (20) can be rewritten as

$$(21) \qquad \sigma(\tau, h) = 2\tau \int_0^2 \frac{a(\alpha) b(\alpha) \mathcal{R}(\alpha)}{h^\alpha} d\alpha \leq 1.$$

It follows from the latter inequality that $h \to 0$ yelds $\tau \to 0$. This, in turn, yields $n = t/\tau \to \infty$ for any finite $t$.

Now we assume that the singular support of $a$ does not contain 2, i.e., $\{2\} \notin \text{singsupp}\, a$ [1].

---

[1] This condition relates only to distributions.

**Theorem 3.1.** *Let* $\mathbf{X}$ *be a random vector with the transition probabilities* $p_k = P(\mathbf{X} = x_k), k \in \mathbb{Z}^N$, *defined in Eqs: (18), (19) and, which satisfy the condition (20) (or, the same, (21)). Then the sequence of random vectors* $\mathbf{S}_n = h\mathbf{X}_1 + \cdots + h\mathbf{X}_n$, *converges as* $n \to \infty$ *in law to the random vector whose probability density function is the fundamental solution of the distributed space fractional order differential equation (10).*

*Proof.* We have to show that the sequence of random vectors $\mathbf{S}_n$ tends to the random vector with pdf $G(t, x)$ in Eq. (15), or the same, the discrete function $y_j(t_n)$ tends to $G(t, x)$ as $n \to \infty$. It is obvious that the Fourier transform of $G(t, x)$ with respect to the variable $x$ is the function $\hat{G}(t, \xi) = e^{t\mathcal{B}(\xi)}$, where $\mathcal{B}(\xi)$ is defined in Eq. (14). Let $\hat{p}(-\xi)$ be the characteristic function corresponding to the discrete function $p_k$, $k \in \mathbb{Z}^N$, that is

$$\hat{p}(-\xi) = \sum_{k \in \mathbb{Z}^N} p_k e^{ik\xi}.$$

It follows from the recursion formula (17) (which exhibits the convolution) and the well known fact that convolution goes over in multiplication by the Fourier transform, the characteristic function of $y_j(t_n)$ can be represented in the form

$$\hat{y}_j(t_n, -\xi) = \hat{p}^n(-\xi).$$

Taking this into account it suffices to show that

(22) $$\hat{p}^n(-h\xi) \to e^{t\mathcal{B}(\xi)}, n \to \infty.$$

The next step of the proof is based on the following simple fact: if a sequence $s_n$ converges to $s$ for $n \to \infty$, then

(23) $$\lim(1 + \frac{s_n}{n})^n = e^s.$$

We have

(24) $$\begin{aligned} \hat{p}^n(-h\xi) &= (1 - \tau \sum_{k \neq 0} \frac{Q_k}{|k|^N}(1 - e^{ik\xi h}))^n \\ &= (1 - \tau \sum_{k \neq 0} \frac{1}{|k|^N} \int_0^2 \frac{a(\alpha)b(\alpha)d\alpha}{|k|^\alpha h^\alpha}(1 - e^{ik\xi h}))^n \\ &= (1 + \frac{t \int_0^2 a(\alpha)\{b(\alpha) \sum \frac{\Delta^2 e^{ik\xi h}}{|kh|^{N+\alpha}} h^N\}d\alpha}{n})^n \end{aligned}$$

It follows from (3) and Corollary 2.5 that

$$b(\alpha) \sum_{k \in Z^N} \frac{\Delta^2 e^{ik\xi h}}{|kh|^{N+\alpha}} h^N$$

tends to $(D_0^\alpha e^{ix\xi})_{|x=0} = -|\xi|^\alpha$ as $h \to 0$ (or, the same, $n \to \infty$) for all $\alpha \in (0, 2]$. Hence

$$s_n = \int_0^2 a(\alpha)\{b(\alpha) \sum \frac{\Delta^2 e^{ik\xi h}}{|kh|^{N+\alpha}} h^N\}d\alpha \to \mathcal{B}(\xi), \ n \to \infty \ (h \to 0).$$

Thus, in accordance with (23) we have

$$\hat{p}^n(-h\xi) \to e^{t\mathcal{B}(\xi)}, n \to \infty. \qquad \square$$

The random walk related to the multiterm fractional diffusion equation can be derived from Theorem 3.1. Assume that $a(\alpha)$ has the form (12) with $0 < \alpha_1 < \cdots < \alpha_M < 2$. So, we again exclude the case $\{2\} \in sing\,supp\,a$.

**Theorem 3.2.** *Let the transition probabilities $p_k = P(\mathbf{X} = x_k), k \in \mathbb{Z}^N$, of the random vector $\mathbf{X}$ be given as follows:*

$$p_k = \begin{cases} 1 - \sum_{j \neq 0} \frac{1}{|j|^N} \sum_{m=1}^{M} \frac{\mu_m a_m b(\alpha_m)}{|j|^{\alpha_m}}, & \text{if } k = 0; \\ \frac{1}{|k|^N} \sum_{m=1}^{M} \frac{\mu_m a_m b(\alpha_m)}{|j|^{\alpha_m}}, & \text{if } k \neq 0, \end{cases}$$

*where $\mu_m = \frac{2\tau}{h^{\alpha_m}}, m = 1, \ldots, M$. Assume,*

$$\sum_{m=1}^{M} a_m b(\alpha_m) \mathcal{R}(\alpha_m) \mu_m \leq 1.$$

*Then the sequence of random vectors $\mathbf{S}_n = h\mathbf{X}_1 + \cdots + h\mathbf{X}_n$, converges as $n \to \infty$ in law to the random vector whose probability density function is the fundamental solution of the multiterm fractional order differential equation (13).*

**Remark 3.3.** The condition $\{2\} \notin sing\,supp\,a$ is required due to singularity of the value $\alpha = 2$ in the definition of $D_0^\alpha$ (see (7)). The particular case $a(\alpha) = \delta(\alpha - 2)$ reduces to the classic heat conduction equation and corresponding random walk is the classic Brownian motion. In more general case of $a(\alpha) = \sum_{l=0}^{m} c_l \delta^{(l)}(\alpha - 2)$ this condition leads to the scaling limit with $\sigma(\tau, h) = h^2 ln\frac{1}{h}$ (see, also [16]).

This work was supported in part by NIH grant P20 GMO67594.

# References

[1] ADLER, R., FELDMAN, R. AND TAQQU, M. (1998). *A Practical Guide to Heavy Tails*. Birkhäuser, Boston. MR1652283
[2] ANDERSON, P. AND MEERSCHAERT, M. M. (1998). Modeling river flows with heavy tails. *Water Resour. Res.* **34** 2271–2280.
[3] ANDRIES, E., STEINBERG, S. AND UMAROV, S. Fractional space-time differential equations: theoretical and numerical aspects, in preparation.
[4] BAGLEY, R. L. AND TORVIC, P. J. (2000). On the existance of the order domain and the solution of distributed order equations I, II. *Int. J. Appl. Math* **2** 865–882, 965–987. MR1760547
[5] CAPUTO, M. (1967). Linear models of dissipation whose Q is almost frequency independent. II. *Geophys. J. R. Astr. Soc.* **13** 529–539.
[6] CAPUTO, M. (2001). Distributed order differential equations modeling dielectric induction and diffusion. *Fract. Calc. Appl. Anal.* **4** 421–442. MR1874477
[7] CHECHKIN, A. V. AND GONCHAR, V. YU. (1999). A model for ordinary Levy motion. *ArXiv:cond-mat/9901064* **v1**, 14 pp.
[8] CHECHKIN, A. V., GORENFLO, R., SOKOLOV, I. M. AND GONCHAR, V. YU. (2003). Distributed order time fractional diffusion equation. *FCAA* **6** 259–279. MR2035651
[9] DIETHELM, K. AND FORD, N. J. (2001). Numerical solution methods for distributed order differential equations. *FCAA* **4** 531–542. MR1874482

[10] DUBINSKII, YU. A. (1991). *Analytic Pseudo-differential Operators and Their Applications.* Kluwer Academic Publishers, Dordrecht. MR1175753
[11] FELLER, W. (1974). *An Introduction to Probability Theory and Its Applications.* John Wiley and Sons, New York–London–Sydney. MR
[12] GILLIS, G. E. AND WEISS, G. H. (1970). Expected number of distinct sites visited by a random walk with an infinite variance. *J. Math. Phys.* **11** 1307–1312. MR0260036
[13] GNEDENKO, B. V. AND KOLMOGOROV, A. N. (1954). *Limit Distributions for Sums of Independent Random Variables.* Addison-Wesley, Reading.
[14] GORENFLO, R., LUCHKO, YU. AND UMAROV, S. (2000). On the Cauchy and multi-point value problems for partial pseudo-differential equations of fractional order. *FCAA* **3** 250–275. MR1788164
[15] GORENFLO, R. AND MAINARDI, F. (1999). Approximation of Lévy-Feller diffusion by random walk. *ZAA* **18** (2) 231–246. MR1701351
[16] GORENFLO, R. AND MAINARDI, M. (2001). Random walk models approximating symmetric space-fractional diffusion processes. In: Elschner, Gohberg and Silbermann (Eds), *Problems in Mathematical Physics* (Siegfried Prössdorf Memorial Volume). Birkhäuser Verlag, Boston-Basel-Berlin, pp. 120–145.
[17] GORENFLO, R., MAINARDI, F., MORETTI, D., PAGNINI, G. AND PARADISI, P. (2002). Discrete random walk models for space-time fractional diffusion. *Chemical Physics.* **284** 521–541.
[18] HÖRMANDER, L. (1983). *The Analysis of Linear Partial Differential Operators: I. Distribution Theory and Fourier Analysis.* Springer-Verlag, Berlin–Heidelberg–New York–Tokyo. MR0717035
[19] JACOB, N. (2001). *Pseudo-differential Operators and Markov Processes.* Vol. **I**: *Fourier Analysis and Semigroups.* Imperial College Press, London. MR1873235
[20] JACOB, N. (2002). *Pseudo-differential Operators and Markov Processes.* Vol. **II**: *Generators and Their Potential Theory.* Imperial College Press, London. MR1917230
[21] JACOB, N. (2005). *Pseudo-differential Operators and Markov Processes.* Vol. **III**: *Markov Processes and Applications.* Imperial College Press, London. MR2158336
[22] LORENZO, C. F. AND HARTLEY, T. T. (2002). Variable order and distributed order fractional operators. *Nonlinear Dynamics* **29** 57–98. MR1926468
[23] MAINARDI, F., LUCHKO, YU. AND PAGNINI, G. (2001). The fundamental solution of the space-time fractional diffusion equation. *FCAA* **4** 153–192. MR
[24] MCCULLOCH, J. (1996). Financial applications of stable distributions. In: *Statistical Methods in Finance: Handbook of Statistics* **14** (G. Madfala and C.R. Rao, Eds). Elsevier, Amsterdam, pp. 393–425. MR1602156
[25] MEERSCHAERT, M. M. AND SCHEFFLER, P.-H. (2001). *Limit Distributions for Sums of Independent Random Vectors. Heavy Tails in Theory and Practice.* John Wiley and Sons, Inc. MR1840531
[26] METZLER, R. AND KLAFTER, J. (2000). The random walk's guide to anomalous diffusion: a fractional dynamics approach. *Physics Reports* **339** 1–77. MR1809268
[27] NIKOL'SKII, S. M. (1975). *Approximation of Functions of Several Variables and Imbedding Theorems.* Springer-Verlag. MR0374877
[28] RUBIN, B. (1996). *Fractional Integrals and Potentials.* Addison Wesley, Longman Ltd. MR1428214
[29] SAMKO, S. G., KILBAS, A. A. AND MARICHEV, O. I. (1993). *Fractional Integrals and Derivatives: Theory and Applications.* Gordon and Breach Science

Publishers, New York and London (Originally published in 1987 in Russian). MR1347689

[30] SAMORODNITSKY, G. AND TAQQU, M. (1994). *Stable Non-Gaussian Random Processes*. Chapman and Hall, New York–London.

[31] SCHMITT, F. G. AND SEURONT, L. (2001). Multifractal random walk in copepod behavior. *Physica* A **301** 375–396.

[32] SOBOLEV, S. L. (1992). *Introduction to the Theory of Cubature Formulas*. Nauka, Moscow (1974) (in Russian). Translated into English: S. L. Sobolev, *Cubature Formulas and Modern Analysis: An Introduction.* Gordon and Breach Science Publishers. MR1248825

[33] UCHAIKIN, V. V. AND ZOLOTAREV, V. M. (1999). *Chance and Stability. Stable Distributions and Their Applications*. VSP, Utrecht. MR1745764

[34] UMAROV, S. (1997, 1998). Non-local boundary value problems for pseudo-differential and differential-operator equations I, II. *Differential Equations* **33, 34** 831–840, 374–381. MR1615099, MR1668214

[35] UMAROV, S. (2003). Multidimensional random walk model approximating fractional diffusion processes. *Docl. Ac. Sci. of Uzbekistan*.

[36] UMAROV, S. AND GORENFLO, R. (2005). On multi-dimensional symmetric random walk models approximating fractional diffusion processes. *FCAA* **8** 73–88.

[37] UMAROV, S. AND GORENFLO, R. (2005). The Cauchy and multipoint problem for distributed order fractional differential equations. *ZAA* **24** 449–466.

[38] ZASLAVSKY, G. (2002). Chaos, fractional kinetics, and anomalous transport. *Physics Reports* **371** 461–580. MR1937584